\documentclass[a4paper,12pt]{article}
\usepackage{amssymb}
\usepackage{amsthm}
\usepackage{amsxtra}
\usepackage{amsfonts}
\usepackage{xr}

\begin{document}

\def\sect{\section}

\newtheorem{thm}{Theorem}[section]
\newtheorem{cor}[thm]{Corollary}
\newtheorem{lem}[thm]{Lemma}
\newtheorem{prop}[thm]{Proposition}
\newtheorem{propconstr}[thm]{Proposition-Construction}

\theoremstyle{definition}
\newtheorem{para}[thm]{}
\newtheorem{ax}[thm]{Axiom}
\newtheorem{conj}[thm]{Conjecture}
\newtheorem{defn}[thm]{Definition}
\newtheorem{notation}[thm]{Notation}
\newtheorem{rem}[thm]{Remarks}
\newtheorem{remark}[thm]{Remark}
\newtheorem{question}[thm]{Question}
\newtheorem{example}[thm]{Example}
\newtheorem{problem}[thm]{Problem}
\newtheorem{excercise}[thm]{Exercise}
\newtheorem{ex}[thm]{Exercise}

\def\Bbb{\mathbb}
\def\cal{\mathcal}
\def\mL{{\mathcal L}}
\def\mC{{\mathcal C}}

\overfullrule=0pt

\def\si{\sigma}
\def\prf{\smallskip\noindent{\it        Proof}. }
\def\call{{\cal L}}
\def\nat{{\Bbb  N}}
\def\la{\langle}
\def\ra{\rangle}
\def\inv{^{-1}}
\def\ld{{\rm    ld}}
\def\trdeg{{tr.deg}}
\def\dim{{\rm   dim}}
\def\th{{\rm    Th}}
\def\rest{{\lower       .25     em      \hbox{$\vert$}}}
\def\ch{{\rm    char}}
\def\zee{{\Bbb  Z}}
\def\conc{^\frown}
\def\acl{acl_\si}
\def\cls{cl_\si}
\def\cals{{\cal S}}
\def\mult{{\rm  Mult}}
\def\calv{{\cal V}}
\def\aut{{\rm   Aut}}
\def\ffi{{\Bbb  F}}
\def\ffiti{\tilde{\Bbb          F}}
\def\degs{deg_\si}
\def\calx{{\cal X}}
\def\gal{{\cal G}al}
\def\cl{{\rm cl}}
\def\loc{{\rm locus}}
\def\calg{{\cal G}}
\def\calq{{\cal Q}}
\def\calr{{\cal R}}
\def\caly{{\cal Y}}
\def\aff{{\Bbb A}}
\def\cali{{\cal I}}
\def\calu{{\cal U}}
\def\epsilon{\varepsilon} 
\def\Uu{{\cal U}}
\def\rat{{\Bbb Q}}
\def\ga{{\Bbb G}_a}
\def\gm{{\Bbb G}_m}
\def\cee{{\Bbb C}}
\def\ree{{\Bbb R}}
\def\frob{{\rm Frob}}
\def\Frob{{\rm Frob}}
\def\fix{{\rm Fix}}
\def\Uu{{\cal U}}
\def\proj{{\Bbb P}}
\def\sym{{\rm Sym}}
 
\def\dcl{{\rm dcl}}
\def\calm{{\mathcal M}}

\font\helpp=cmsy5
\def\semdp
{\hbox{$\times\kern-.23em\lower-.1em\hbox{\helpp\char'152}$}\,}

\def\dnfo{\,\raise.2em\hbox{$\,\mathrel|\kern-.9em\lower.35em\hbox{$\smile$}
$}}
\def\dnf#1{\lower1em\hbox{$\buildrel\dnfo\over{\scriptstyle #1}$}}
\def\dfo{\;\raise.2em\hbox{$\mathrel|\kern-.9em\lower.35em\hbox{$\smile$}
\kern-.7em\hbox{\char'57}$}\;}
\def\df#1{\lower1em\hbox{$\buildrel\dfo\over{\scriptstyle #1}$}}        
\def\stab{{\rm Stab}}
\def\qfcb{\hbox{qf-Cb}}
\def\perf{^{\rm perf}}
\def\sipm{\si^{\pm 1}}

\newcommand{\nc}{\newcommand}
\nc{\renc}{\renewcommand}
\nc{\ssec}{\subsection}
\nc{\sssec}{\subsubsection}
\nc{\on}{\operatorname}

\nc\ol{\overline}
\nc\wt{\widetilde}
\nc\wh{\widehat}
\nc\tboxtimes{\wt{\boxtimes}}

\emergencystretch=2cm

\nc{\Aa}{{\mathbb{A}}}
 \nc{\Gg}{{\mathbb{G}}}  \def\gg{\mathbb{g}}
\nc{\Hh}{{\mathbb{H}}}
 \nc{\Nn}{{\mathbb{N}}}
\nc{\Pp}{{\mathbb{P}}}
\nc{\Rr}{{\mathbb{R}}}
\nc{\BV}{{\mathbb{V}}}
\nc{\BW}{{\mathbb{W}}}
\nc{\Zz}{{\mathbb{Z}}}
\nc{\Qq}{{\mathbb{Q}}}
\nc{\Ss}{{\mathbb{S}}}
\nc{\Cc}{{\mathbb{C}}}

\nc{\CA}{{\mathcal{A}}}
\nc{\CB}{{\mathcal{B}}}

\nc{\CE}{{\mathcal{E}}}
\nc{\CF}{{\mathcal{F}}}
\nc{\CG}{{\mathcal{G}}}
\nc{\CL}{{\mathcal{L}}}
\nc{\CC}{{\mathcal{C}}}
\nc{\CM}{{\mathcal{M}}}
\def\Mm{\CM}
\nc{\CN}{{\mathcal{N}}}
\nc{\Oo}{{\mathcal{O}}}
\nc{\CP}{{\mathcal{P}}}
\nc{\CQ}{{\mathcal{Q}}}
\nc{\CR}{{\mathcal{R}}}
\nc{\CS}{{\mathcal{S}}}
\nc{\CT}{{\mathcal{T}}}
\nc{\CU}{{\mathcal{U}}}
\nc{\CV}{{\mathcal{V}}}
\nc{\CK}{{\mathcal{K}}}
\nc{\CW}{{\mathcal{W}}}
\nc{\CZ}{{\mathcal{Z}}}

\nc{\cM}{{\check{\mathcal M}}{}}
\nc{\csM}{{\check{\mathcal A}}{}}
\nc{\oM}{{\overset{\circ}{\mathcal M}}{}}
\nc{\obM}{{\overset{\circ}{\mathbf M}}{}}
\nc{\oCA}{{\overset{\circ}{\mathcal A}}{}}
\nc{\obA}{{\overset{\circ}{\mathbf A}}{}}
\nc{\ooM}{{\overset{\circ}{M}}{}}
\nc{\osM}{{\overset{\circ}{\mathsf M}}{}}
\nc{\vM}{{\overset{\bullet}{\mathcal M}}{}}
\nc{\nM}{{\underset{\bullet}{\mathcal M}}{}}
\nc{\oD}{{\overset{\circ}{\mathcal D}}{}}
\nc{\obD}{{\overset{\circ}{\mathbf D}}{}}
\nc{\oA}{{\overset{\circ}{\mathbb A}}{}}
\nc{\op}{{\overset{\bullet}{\mathbf p}}{}}
\nc{\cp}{{\overset{\circ}{\mathbf p}}{}}
\nc{\oU}{{\overset{\bullet}{\mathcal U}}{}}
\nc{\oZ}{{\overset{\circ}{\mathcal Z}}{}}
\nc{\ofZ}{{\overset{\circ}{\mathfrak Z}}{}}
\nc{\oF}{{\overset{\circ}{\fF}}}

\nc{\fa}{{\mathfrak{a}}}
\nc{\fb}{{\mathfrak{b}}}
\nc{\fg}{{\mathfrak{g}}}
\nc{\fgl}{{\mathfrak{gl}}}
\nc{\fh}{{\mathfrak{h}}}
\nc{\fj}{{\mathfrak{j}}}
\nc{\fm}{{\mathfrak{m}}}
\nc{\fn}{{\mathfrak{n}}}
\nc{\fu}{{\mathfrak{u}}}
\nc{\fp}{{\mathfrak{p}}}
\nc{\fr}{{\mathfrak{r}}}
\nc{\fs}{{\mathfrak{s}}}
\nc{\fsl}{{\mathfrak{sl}}}
\nc{\hsl}{{\widehat{\mathfrak{sl}}}}
\nc{\hgl}{{\widehat{\mathfrak{gl}}}}
\nc{\hg}{{\widehat{\mathfrak{g}}}}
\nc{\chg}{{\widehat{\mathfrak{g}}}{}^\vee}
\nc{\hn}{{\widehat{\mathfrak{n}}}}
\nc{\chn}{{\widehat{\mathfrak{n}}}{}^\vee}

\nc{\fA}{{\mathfrak{A}}}
\nc{\fB}{{\mathfrak{B}}}
\nc{\fD}{{\mathfrak{D}}}
\nc{\fE}{{\mathfrak{E}}}
\nc{\fF}{{\mathfrak{F}}}
\nc{\fG}{{\mathfrak{G}}}
\nc{\fK}{{\mathfrak{K}}}
\nc{\fL}{{\mathfrak{L}}}
\nc{\fM}{{\mathfrak{M}}}
\nc{\fN}{{\mathfrak{N}}}
\nc{\fP}{{\mathfrak{P}}}
\nc{\fU}{{\mathfrak{U}}}
\nc{\fV}{{\mathfrak{V}}}
\nc{\fZ}{{\mathfrak{Z}}}

\nc{\bb}{{\mathbf{b}}}
\nc{\bc}{{\mathbf{c}}}
\nc{\bd}{{\mathbf{d}}}
\nc{\be}{{\mathbf{e}}}
\nc{\bj}{{\mathbf{j}}}
\nc{\bn}{{\mathbf{n}}}
\nc{\bp}{{\mathbf{p}}}
\nc{\bq}{{\mathbf{q}}}
\nc{\bF}{{\mathbf{F}}}
\nc{\bu}{{\mathbf{u}}}
\nc{\bv}{{\mathbf{v}}}
\nc{\bx}{{\mathbf{x}}}
\nc{\bs}{{\mathbf{s}}}
\nc{\by}{{\mathbf{y}}}
\nc{\bw}{{\mathbf{w}}}
\nc{\bA}{{\mathbf{A}}}
\nc{\bK}{{\mathbf{K}}}
\nc{\bI}{{\mathbf{I}}}
\nc{\bB}{{\mathbf{B}}}
\nc{\bG}{{\mathbf{G}}}
\nc{\bC}{{\mathbf{C}}}
\nc{\bD}{{\mathbf{D}}}
\nc{\bP}{{\mathbf{P}}}
\nc{\bH}{{\mathbf{H}}}
\nc{\bM}{{\mathbf{M}}}
\nc{\bN}{{\mathbf{N}}}
\nc{\bV}{{\mathbf{V}}}
\nc{\bU}{{\mathbf{U}}}
\nc{\bL}{{\mathbf{L}}}
\nc{\bT}{{\mathbf{T}}}
\nc{\bW}{{\mathbf{W}}}
\nc{\bX}{{\mathbf{X}}}
\nc{\bY}{{\mathbf{Y}}}
\nc{\bZ}{{\mathbf{Z}}}
\nc{\bS}{{\mathbf{S}}}
\nc{\bQ}{{\mathbf{Q}}}

\nc{\sA}{{\mathsf{A}}}
\nc{\sB}{{\mathsf{B}}}
\nc{\sC}{{\mathsf{C}}}
\nc{\sD}{{\mathsf{D}}}
\nc{\sF}{{\mathsf{F}}}
\nc{\sG}{{\mathsf{G}}}
\nc{\sK}{{\mathsf{K}}}
\nc{\sM}{{\mathsf{M}}}
\nc{\sO}{{\mathsf{O}}}
\nc{\sQ}{{\mathsf{Q}}}
\nc{\sP}{{\mathsf{P}}}
\nc{\sZ}{{\mathsf{Z}}}
\nc{\sfp}{{\mathsf{p}}}
\nc{\sr}{{\mathsf{r}}}
\nc{\sg}{{\mathsf{g}}}
\nc{\sff}{{\mathsf{f}}}
\nc{\sfb}{{\mathsf{b}}}
\nc{\sfc}{{\mathsf{c}}}
\nc{\sd}{{\ltimes}}

\nc{\tA}{{\widetilde{\mathbf{A}}}}
\nc{\tB}{{\widetilde{\mathcal{B}}}}
\nc{\tg}{{\widetilde{\mathfrak{g}}}}
\nc{\tG}{{\widetilde{G}}}
\nc{\TM}{{\widetilde{\mathbb{M}}}{}}
\nc{\tO}{{\widetilde{\mathsf{O}}}{}}
\nc{\tU}{{\widetilde{\mathfrak{U}}}{}}
\nc{\TZ}{{\tilde{Z}}}
\nc{\tx}{{\tilde{x}}}
\nc{\tq}{{\tilde{q}}}

\nc{\tfP}{{\widetilde{\mathfrak{P}}}{}}
\nc{\tz}{{\tilde{\zeta}}}
\nc{\tmu}{{\tilde{\mu}}}

\newcommand\Hom{{\rm Hom}}

  \nc{\Ob}{{\mathop{\operatorname{\rm Ob}}}}
  \nc{\Sym}{{\mathop{\operatorname{\rm Sym}}}}
   \nc{\Aut}{{\mathop{\operatorname{\rm Aut}}}}
 \nc{\Spec}{{\mathop{\operatorname{\rm Spec}}}}
  \nc{\spec}{{\mathop{\operatorname{\rm Spec}}}}
\nc{\Ker}{{\mathop{\operatorname{\rm Ker}}}}
 \nc{\dom}{{\mathop{\operatorname{\rm dom}}}}
\nc{\End}{{\mathop{\operatorname{\rm End}}}}
 \nc{\GL}{{\mathop{\operatorname{\rm GL}}}}
 \nc{\Id}{{\mathop{\operatorname{\rm Id}}}}
 \nc{\rk}{{\mathop{\operatorname{\rm rk}}}} 
 \nc{\length}{{\mathop{\operatorname{\rm length}}}}
\nc{\supp}{{\mathop{\operatorname{\rm supp}}}}
\nc{\val}{{\rm val}}
\nc{\res}{{\mathop{\operatorname{\rm res}}}}
\def\ind#1#2{ {#1} {\downarrow} {#2} }  
\def\Ind#1#2#3{{#1} {\downarrow}_{#3} {#2} }  
\def\domn#1#2#3{ (#1 \leq_{dom_{#3}} #2 ) }
\def\tensor{{\otimes}}
\def\meet{\cap}
\def\union{\cup}
\def\Union{\bigcup}
\def\si{\sigma}
\def\g{\gamma}
\def\G{\Gamma}
\def\Sum{\Sigma}
\def\<{\begin}
 \def\>{\end}
\def\m{\setminus}

 \def\AFD{(FD)\,}
\def\Aom{(FD$_{\omega}$)\,}
\nc{\seq}[1]{\stackrel{#1}{\sim}}
\def\dd#1{\frac{\partial}{\partial X_{#1}}}   \def\ddt{\frac{\partial}{\partial t}}
\def\invv#1{{{#1}^{-1}}}
\def\inv{^{-1}}
\def\claim#1{{\noindent \bf Claim #1.\ }}
\def\Claim{{\noindent \bf Claim.\ }}
\def\beq#1{{\begin{equation} \label{#1}}  }
\def\normal{\trianglelefteq}
\def\Uu{\mathbb U}
\def\iso{\simeq}

\def\pv{\hfill $\Box$}
\def\eprf{\hfill $\Box$}
 
\def\acl{\mathop{\rm acl}\nolimits}
 \def\dcl{\mathop{\rm dcl}\nolimits}
\def\OK{{\mathcal O}_K}
\def\liminv{\underset{\longleftarrow}{lim}\,} 
\def\liminvi{\underset{\longleftarrow i}{lim}\,} 

\def\lbl#1{    \label{#1}  }
\def\a{\alpha}
\def\revo{\backslash}
\def\ba{\bar{a}}
\def\bba{{\mathbf a}}
\def\k{{\rm k}}
 
\def\tlt{\times \ldots \times}

  \def\Rr{\ree}

\def\frac#1#2{{#1\over #2}}

\def\lam{\lambda}
\def\e{\epsilon}

\def\vlabel{\label}

\def\udirem{\smallskip\noindent {\bf Remark}. }

\title{Difference fields and descent in algebraic dynamics, II}

\author{Zo\'e Chatzidakis\thanks{partially supported by MRTN-CT-2004-512234 and by ANR-06-BLAN-0183}{ (CNRS - Universit\'e Paris 7)}
\and
Ehud Hrushovski\thanks{thanks to the Israel Science Foundation for
  support, and to the ANR-06-BLAN-0183}\  (The Hebrew University)}
\date{}
\maketitle
 \def\hK{{\widehat{K}}}
\def\hL{{\widehat{L}}}

\section*{Introduction}
This second part of the paper strengthens  the descent theory described
in the first part \cite{dyn1}  to
rational maps and arbitrary base fields.  In particular this is needed in order to obtain the
``dynamical Northcott'' Theorem 1.11 of Part I in sharp form.  As a special case, for   $\Pp^1$
over an algebraically closed field, we recover the theorem of \cite{baker}.
 
 Stated model-theoretically, our
main result is the following:

\bigskip\noindent
{\bf Theorem \ref{theo1a}} {\em Let $(\calu,\si)$ be an inversive  difference field, and $K_1\subset K_2$ be subfields of the fixed field
  $\fix(\si)=\{c\in\calu\mid \si(c)=c\}$, with $K_2/K_1$ regular. Let
  $a,b$ be tuples in $\calu$ 
  such that $SU(a/K_2)<\infty$ and: \begin{itemize}
\item[\rm (a)] $a$ belongs to the  difference field $K_2(b)_{\si}$
  generated by $b$ over $K_2$; 
\item [\rm (b)]$K_1(b)_{\si}$ is linearly disjoint from
  $K_2$ over $K_1$; 
\item [\rm (c)]$tp(a/K_2)$ is hereditarily orthogonal to
  $\fix(\si)$. 
\end{itemize} 

Then there is a tuple $c $ in $K_1(b)\perf_\si$ such that $a\in K_2(c)_\si$ and $c$ is
purely inseparable over $K_2(a)_\si$. If $tr.deg(K_2(a)_\si/K_2)=1$ and
$K_1$ is perfect, 
then one can choose $c$ so that  $K_1(b)_\si\cap K_2(a)_\si=K_1(c)_\si$.}

\bigskip
This result is in fact an easy corollary of a similar statement for {\it
  inversive} difference fields, see Proposition \ref{theo1}. Stated in the language of the first part \cite{dyn1}, the orthogonality
condition simply means {\em fixed-field-free}, but see also section
\ref{prelim} for a more model-theoretic definition. 
This theorem implies a result which can be stated in geometric
terms, given the language of algebraic dynamics (cf. the introduction of 
\cite{dyn1}): 

\bigskip\noindent
{\bf Theorem 
\ref{theo2}} {\em  Let $K_1\subset K_2$ be fields, with $K_2/K_1$
regular,  and let
  $(V_2,\phi_2)\in AD_{K_2}$. Assume that $(V_2,\phi_2)$ is
  primitive and that
  $\deg(\phi_2)>1$. Assume furthermore that   for some $n\geq 1$,  $(V_2,\phi_2^n)$ is dominated
  (in $AD_{K_2}$)
  by some object of $AD_{K_1}$.
\begin{enumerate}

\item There is some variety $V_3$ defined over $K_1$, and a
  dominant constructible map $\phi_3:V_3\to V_3$ also defined over $K_1$, a constructible
  isomorphism $h:(V_2,\phi_2)\to (V_3,\phi_3)$. 
\item Assume that the characteristic is $0$, or that $K_1$ is perfect
  and $\dim(V_2)=1$. Then $(V_2,\phi_2)$ is rationally isotrivial, i.e., there is
  some $(V_3,\phi_3)\in
AD_{K_1}$ which is isomorphic to $(V_2,\phi_2)$ (in $AD_{K_2}$).
\end{enumerate}
}
Note that in the general case (1) above, we  do not know that
$(V_3,\phi_3)$ can be taken in $AD_{K_1}$: we only know that $\phi_3$ is
constructible. The statement of Theorem \ref{theo2} was devised to contain no model-theoretical notions, but
the proof will really 
use only that $(V_2,\phi_2)$ is fixed-field-free; in fact even this
assumption is not needed here, as will be seen in \cite{dyn3} (see point
1 below).

Before describing the organisation of the paper and explaining some of the tools arising in the proof, here are
some of the highlights of the third paper \cite{dyn3}. Indeed,  Proposition \ref{theo1}
and Theorem \ref{theo2} suggest 
 several directions of research.

\smallskip
{\bf 1. } The assumption on $\deg(\phi_2)$ made in  Theorem
\ref{theo2} is unnecessary:

\smallskip\noindent
{\bf Theorem}. {\em Let $K_1=K_1^{alg}\subset K_2$ be fields,
  $(V_2,\phi_2)\in AD_{K_2}$ with $V_2$  
absolutely irreducible. Assume that $(V_2,\phi_2)$ is dominated by an
object of $AD_{K_1}$. Then $(V_2,\phi_2)$ dominates some $(V_3,\phi_3)
\in AD_{K_1}$ (with $\dim(V_3)>0$). If $(V_2,\phi_2)$ is  primitive, it is (constructibly) isotrivial.}

Note the extra assumption of $K_1$ being algebraically closed. The case not covered by the proof of Proposition \ref{theo1} is when
$tp(a/K_1)$ is almost internal to $\fix(\si)$. We  use the Galois theory of
difference equations (a standard model-theoretic tool, which has to be
slightly adapted to our context; see also the
end of the introduction of the first part \cite{dyn1}). Assuming $K_2$
algebraically closed, we show that $a$ is equi-algebraic over 
$K_2$ to an element $b\in K_2(a)_{\si^{\pm 1}}$ in a translation variety, more
precisely 
satisfying an equation $x\in G \land \si(x)=x+g$, where $G$ is a
(simple) commutative algebraic group and $g\in G(K_2)$. From this the
result follows, using the Galois correspondence. 

\smallskip
{\bf 2.}   The
hypothesis of hereditary orthogonality in Proposition \ref{theo1} cannot 
be removed; see  \ref{counterex}. However,
using a model-theoretic result on canonical bases, one can show that
there is some $c\in K_2(a)_{\si^{\pm 1}}^{alg}$ such that $tp(a/K_2(c)_{\si^{\pm 1}})$ is almost
internal to $\fix(\si)$, and $K_1(c)_{\si^{\pm 1}}$ is linearly disjoint from
$K_2$ over $K_1$. In particular, if $a$ is the generic of some
$(V,\phi)\in AD_{K_2}$, then $c$ is the generic of some $(W,\psi)\in
AD_{K_1}$ with $\deg(\psi)=\deg(\phi)$. It is unlikely that such a $c$
can always be found in $K_2(a)_{\si^{\pm 1}}$, but we do not know of  a
counterexample.  

\smallskip
{\bf 3.} All results mentioned above generalize to the case where $K_1$
is not necessarily a subfield of $K_2$: for instance, in \ref{theo2},
the object $(V_3,\phi_3)$ will be defined over $K_0=K_1^{alg}\cap
K_2^{alg}$.  If $K_0\subset K_2$, then the
(iso)morphism $(V_2,\phi_2)\to (V_3,\phi_3)$ is also defined over $K_2$.


\bigskip
The first section of the paper sets up the notation and recalls some of the classical results
from stability theory and the model theory of existentially closed difference fields. In
section 2, we prove some more technical results on difference fields,
and in particular investigate quantifier-free types in reducts
$(\calu,\si^n)$ of our existentially closed difference field $\calu$. Our two main
results in this section are Proposition \ref{propo4} and Theorem
\ref{theo3}, and 
are we believe of independent interest. Finally, section 3 gives the
proofs of the two main results.

\section{\bf Preliminaries on difference fields}\vlabel{prelim}
\para{\bf Notation}. A {\it difference field} is a
field with a distinguished endomorphism $\si$. We denote by  $\call$ the
language of rings $\{+,-,\cdot, 0,1\}$ and by $\call_\si$ the language
$\call\cup\{\si\}$, where $\si$ is a unary function symbol. Difference
fields are then $\call_\si$-structures. Recall that a difference field
is {\em inversive} if the endomorphism $\si$ is also surjective.

Let $K\subset L$ be difference fields, $E$ a field. If $a$ is a tuple in
$L$, then $K(a)_\si$ will denote the  difference subfield of $L$
generated by $a$ over $K$:  $K(a)_\si=K(\si^i(a)\mid i\in
\nat)$; if $K$ is { inversive}, then 
$K(a)_{\si^{\pm 1}}$ will denote the {inversive} difference subfield of $L$
generated by 
$a$: $K(a)_{\si^{\pm 1}}=K(\si^i(a)\mid i\in\zee)$. We denote by
$E^{sep}$ the separable closure of the field $E$,  by $E^{alg}$ its
algebraic closure, and by $E\perf$ the perfect hull $E^{1/p^\infty}$
over $E$ if ${\rm char}(E)=p>0$, the field $E$ if ${\rm char}(E)=0$. 
If $F$ is a
(finite) algebraic extension of $E$, then 
$[F:E]_s$ will denote the separable degree of $F$ over $E$, and
$[F:E]_i$ its inseparable degree. If the characteristic is $p>0$, then $\frob$
will denote the Frobenius map $x\mapsto x^p$.

\para\vlabel{acfa0}{\bf Existentially closed difference fields, the theory ACFA}. 
The  natural setting for studying difference equations is within
existentially closed 
difference fields.
A  difference field $(\calu,\si)$ is {\it existentially closed}  if any
finite set of difference 
equations (over $\calu$)  with a solution in some  difference field
extension, already has a
solution in $\calu$.   

 The
class of existentially closed difference fields is axiomatisable, and
its theory is 
denoted by ACFA. All completions of ACFA are supersimple. We quickly
recall some of the main properties of ACFA. For details, see
\cite{CH}. Let $\calu$ be an existentially closed difference field, $K$
an inversive
difference subfield of $\calu$, and $a$ a finite tuple of elements of $\calu$.

\begin{enumerate}
\item The {\it completions} of ACFA are obtained by specifying the characteristic
and describing the isomorphism type of the algebraic closure of the
prime field.  
\item
 If $A\subset
\calu$, then 
$\acl(A)$, the model-theoretic algebraic closure of $A$, coincides with 
the 
smallest algebraically closed inversive difference field containing $A$. Thus, for
instance, if  $a$ is a tuple
in $\calu$, then 
$\acl(Ka)=K(a)_{\si^{\pm 1}}^{alg}$.  
\item The
model-theoretic definable closure of $A$ (denoted by $\dcl(A)$) contains
the inversive difference field
generated by $A$, but is usually much larger, even when the characteristic is
$0$. 
\item Independence is defined 
using independence in the theory ACF of algebraically closed fields: $A$ and $B$ are
{\it independent} over $C$ ($A\dnfo_CB$) if and only if $\acl(CA)$ and
$\acl(CB)$ are linearly disjoint over $\acl(C)$. 

\item The
{\it quantifier-free type of $a$ over $K$} describes the isomorphism type over
$K$ of $K(a)_{\si^{\pm 1}}$, i.e.: $qftp(a/K)=qftp(b/K)$ if and only if there is
a $K$-isomorphism of difference fields $K(a)_{\si^{\pm 1}}\to K(b)_{\si^{\pm 1}}$ which
sends $a$ to $b$ (if and only if there is a $K$-isomorphism of
difference fields $K(a)_\si\to K(b)_\si$ which sends $a$ to $b$). 
\item The {\it
type 
of $a$ over $K$} describes the isomorphism type over $K$ of
the algebraic closure of $Ka$: $tp(a/K)=tp(b/K)$ if and only if
there is a $K$-isomorphism of difference fields $\acl(Ka)\to \acl(Kb)$
which sends $a$ to 
$b$. 
\item
If $a$ is an element of $\calu$, then either the elements $\si^i(a)$,
$i\in\zee$,  are
algebraically independent over $K$, in which case $a$ is {\it
transformally transcendental over $K$} and $SU(a/K)=\omega$; or the transcendence degree of
$K(a)_{\si^{\pm 1}}$ over $K$ is finite, i.e., $a$ is {\it transformally algebraic
over $K$}, and then $SU(a/K)<\omega$. Observe that if $a$ is a tuple of
elements which are transformally algebraic over $K$, then
$SU(a/K)<\omega$, and for some $m$, $K(a)_{\si^{\pm 1}}\subset
K(a,\si(a),\ldots,\si^m(a))^{alg}=\acl(Ka)$. In that case, we will often replace $a$ by
$(a,\si(a),\ldots,\si^m(a))$ and assume that $\si(a)\in K(a)^{alg}$.
\item If $K=K^{alg}$ and $L$ is a difference field containing $K$, there
  is a $K$-embedding of $L$ into some elementary extension of $\calu$. 
\end{enumerate}

For more properties of   difference equations and of
the theory ACFA, we
refer to \cite{CH} and \cite{CHP}.

\para {\bf Conventions}. Unless otherwise mentioned, all our difference
fields will be {\bf inversive}. This is in contrast with the first part,
where  the point of view was more geometric and thus it
  was convenient that difference fields arising from algebraic dynamics
  are finitely generated as fields;  in this second part of the paper,
  the objects we are considering are difference fields and their
  algebraic closures, and it is more convenient for our purposes to have
  the smaller field be inversive. Observe for instance
  that if $\alpha$ is algebraic over $K$ and $\si(K)= K$, then
  $[K(\alpha):K]=[K(\si(\alpha)):K]$;  this is not necessarily the
  case if $\si(K)\neq K$.
This is really just a matter of convenience and does not affect the
generality of our results: every difference field has a unique (up to
isomorphism) inversive closure, see Cohn's book \cite{[Co]}, 2.5.II.

Unless otherwise mentioned, the letters $a$, $b$,  $x$, $y$, \dots\ will
denote finite tuples of elements or
variables, and we will abusively write e.g. $a\in K$ instead of $a$ is a
tuple of elements of $K$. 

\smallskip
In the next few lemmas, we fix a difference
subfield $K$ of a {\it  sufficiently saturated} existentially closed difference field
$\calu$.  By sufficiently saturated, we mean that it is
$\kappa$-saturated for some cardinal $\kappa$ greater than all
cardinalities of fields we consider, so
that for instance, in item (8) above the $K$-embedding of $L$ can be
taken  into $\calu$.

\para\vlabel{lem2}{\bf Fixed fields}. A {\it fixed field} is a subfield
of $\calu$ which is 
defined by the equation $\si^n(x)=x$, or by an equation of the form
$\tau(x)=x$ where $\tau=\si^n\frob^m$, $n>0$, $m\in\zee$, if
$char(\calu)=p>0$. The 
fixed field defined by $\tau(x)=x$ is denoted by $\fix(\tau)$; it has
SU-rank $1$ if $n=1$, or if 
$m\neq 0$ and $n,m$ are relatively prime (see 7.1 in \cite{CHP}). These
conditions are clearly necessary, since $\fix(\tau^\ell)$ is an
$\ell$-dimensional $\fix(\tau)$-vector space.

\smallskip\noindent
{\bf Fact} (see the proof of 3.7(3) in \cite{CH}). Let $\tau=\si^n\Frob^m$, and let
$a$ be a tuple in $\fix(\tau)$. Then $K$ and $\fix(\tau)$ are linearly
disjoint over their intersection, and therefore the field of definition of the
algebraic locus of $a$ over $K$ is contained in 
$\fix(\tau)\cap K\perf$. In particular, if $a\in K^{alg}$ and $b$ is the tuple
encoding the set of field conjugates  of $a$ over $K$, then
$b\in \fix(\tau)$. 

\para\vlabel{onebased}{\bf One-basedness}.   Recall that $\calu$
eliminates imaginaries, see \cite{CH}.  
A subset $S$ of $\calu^n$ which is 
invariant under $\aut(\calu/K)$ is {\it one-based (over $K$)} if for any
$K\subset L$ and tuple $a$ of
elements of $S$, $a$ and $L$ are
independent over $\acl(Ka)\cap \acl(L)$. A partial type over $K$ is one-based if
the set of its realisations is one-based (thus any extension of a
one-based type is one-based).  In the presence of a finite dimension theory this is
equivalent to the dimension inequality of Theorem \ref{mod} in \cite{dyn1}.  
 The terms {\em modular}
(used in  \cite{CH}, \cite{CHP}) and {\em one-based} are thus synonymous; we will 
use the former when referring to   $(V,\phi)$ itself, the latter for the solution set
of $\si(x)=\phi(x)$.  
A finite union of one-based sets is
one-based, and the same is true for uniformly definable unions: if $tp(a/K)$
and $tp(b/\acl(Ka))$ are one-based, then so is $tp(a,b/K)$ (see \cite{CH}, \cite{h-emm} or in greater generality \cite{[W]}).   Types of
infinite SU-rank are not one-based, nor are 
(non-algebraic) types which are realised in a fixed field. One  important property
is the following (see e.g. Lemma 3.3 in \cite{CH}):

\smallskip\noindent
{\bf Fact}. Assume that $tp(a/K)$ is one-based, and
$L$ and $K$ are 
independent over $K_0\subset K,L$. Then $\acl(K_0a)$ and $\acl(L)$ are
independent over their intersection. 

\para\vlabel{dich}{\bf The dichotomy}. The main result of \cite{CHP} asserts that if $a$ is a tuple in $\calu$ with
$SU(a/K)=1$, then $tp(a/K)$ is not one-based if and only if $tp(a/K)$ is
non-orthogonal to a fixed field, i.e., there is some algebraically
closed difference field $L$ containing $K$ and linearly disjoint from
$K^{alg}(a)_{\si^{\pm 1}}$ over $K^{alg}$, such that $L(a)_{\si^{\pm 1}}$ contains a tuple
$b\notin L$, $b$ in some fixed field (non-orthogonality gives $b\in
\acl(La)$, use Fact \ref{lem2} to get $b\in L(a)_{\si^{\pm 1}}$). 

\para\vlabel{int1}{\bf Two definitions of internality}. Let $a$ be a
tuple in $\calu$, $a\notin K^{alg}$, and $\pi$ a set of partial types  (over
various parameter sets) which is stable under $\Aut(\calu/K)$.

\begin{enumerate}
\item{We say that $tp(a/K)$ is {\it qf-internal to $\pi$} 
if there
is some $L=\acl(L)$ containing $K$ and independent
from $a$ over $K$, and a tuple $b$ of realisations of types in $\pi$
with base contained in $L$ and  such that
$a\in L(b)_{\si^{\pm 1}}$. (This notion is stronger than the usual notion of
internality, which only requires $a\in \dcl(Lb)$).} 
\item{We say that $tp(a/K)$ is {\it almost internal to $\pi$} if 
there
is some $L=\acl(L)$ containing $K$ and independent
from $a$ over $K$, and a tuple $b$ of realisations of types in $\pi$ with base contained in $L$  such that
$a\in \acl(Lb)$.}
\end{enumerate}

By abuse of language, we also speak of qf-internality (or almost internality) to
$\Pi$, where $\Pi$ is the set of realisations of types in $\pi$.
In
practice, the set $\pi$ will be a union of some of the following sets:
 non-algebraic $1$-types containing $\tau(x)=x$ for some
$\tau=\si^n\frob^m$; all one-based types of SU-rank $1$. 
 
\para\vlabel{int15}{\bf Internality to a fixed field}. Assume that $\Pi$ is the fixed field $\fix(\tau)$. Using Fact
\ref{lem2}, one easily deduces 
\begin{itemize}
\item[(1)]{$tp(a/K)$ is qf-internal to $\fix(\tau)$ if and only if for
some $L$ independent from $a$ over $K$, $L(a)_{\si^{\pm 1}}=L(b)_{\si^{\pm 1}}$ for some
$b\in \fix(\tau)$.}
\item[(2)]{$tp(a/K)$ is almost internal to $\fix(\tau)$ if and only if
for some $L$ independent from $a$ over $K$ and for some tuple $b$ in
$\fix(\tau)\cap L(a)_{\si^{\pm 1}}$, $a\in \acl(Lb)$.} \end{itemize}

\smallskip\noindent
{\bf Remarks}. \begin{enumerate}\item 
Observe that being qf-internal to $\tau(x)=x$ or to
$\tau^k(x)=x$ is the same thing, because $\fix(\tau^k)$ is a
$k$-dimensional vector space over $\fix(\tau)$ and therefore
$\fix(\tau^k)=\fix(\tau)(b)$ if $b$ is any $\fix(\tau)$-basis of
$\fix(\tau^k)$.
\item For algebraic dynamics, this notion is called {\em field-internal}
  in \cite{dyn1}, and {\em fixed-field-internal} if $\tau=\si$.
\end{enumerate}


\para\vlabel{ana}{\bf Analyses}. Let $K=\acl(K)\subset \calu$, $a$ a tuple in
$\calu$, with $SU(a/K)$ finite. A tuple (of tuples) $(a_1,\ldots, a_n)$
is a {\it semi-minimal analysis} of $a$ over $K$ (or of $tp(a/K)$) iff
$\acl(Ka_1,\ldots,a_n)=\acl(Ka)$, and for every $i$,
$tp(a_i/\acl(Ka_1,\ldots,a_{i-1}))$ is (almost- or qf-) internal to the set of
conjugates of a type of
SU-rank $1$.  

By general properties of supersimple theories, every finite SU-rank
type has a semi-minimal analysis. If it is one-based, then one can find
an analysis in which all types $tp(a_i/\acl(Ka_1,\ldots,a_{i-1}))$ have
SU-rank $1$.

The properties of one-based types and the dichotomy (see \ref{onebased}
and \ref{dich}) then yield: 

\smallskip\noindent
{\bf Fact}. A type of finite SU-rank is one-based if and only if it has
a semi-minimal analysis in which all types are one-based types, if and
only if any of its extensions is orthogonal to all fixed fields, if
and only if any of its semi-minimal analyses  only involves types orthogonal to
all fixed fields. 

Recall that a type is {\em hereditarily orthogonal} to  a set $\pi$
of types, if all its extension are orthogonal to all members of
$\pi$. Thus another way of rephrasing the previous fact is: a type is
one-based if and only if it is hereditarily orthogonal to all (types
realised in) fixed fields. In the context of algebraic dynamics, this
corresponds to {\em field-free}, see \cite{dyn1}.

\bigskip

\para{\bf The limit degree and the inverse limit degree}.  These are  numerical invariants that will be helpful in proving
closure properties of $AD_K$ within the category of difference varieties
over $K$. Indeed, an object $(V,\phi)\in AD_K$ will correspond to difference field
extensions with limit degree
$1$ and inverse limit degree $\deg(\phi)$.

\smallskip
\noindent{\bf Definition}. Let $a$ be a tuple in $\calu$, and assume that
$\si(a)\in K(a)^{alg}$. 
The {\it limit degree of $a$ over $K$} (or of $K(a)_{\si^{\pm 1}}$ over
$K$) is $$ld(a/K)=\lim_{k\rightarrow
\infty}[K(a,\si(a),\ldots,\si^{k+1}(a)):K(a,\si(a),\ldots,\si^k(a))],$$
and the {\it inverse limit degree of $a$ over $K$} is
$$ild(a/K)=\lim_{k\rightarrow
\infty}[K(a,\si\inv(a),\ldots,\si^{-(k+1)}(a)):K(a,\si\inv(a),\ldots,\si^{-k}(a))].$$

The limit and inverse limit degrees  are invariants of the
extension $K(a)_{\si^{\pm 1}}/K$, see \cite{[Co]}, section 5.16. If $[K(a,\si(a)):K(a)]=ld(a/K)$, then
the fields $K(\si^i(a)\mid i\leq 0)$ and $K(\si^i(a)\mid i\geq 0)$ are
linearly disjoint over $K(a)$, and
$ild(a/K)=[K(a,\si(a)):K(\si(a))]$. Another important property is that
these degrees are multiplicative in tower. 

\begin{lem}\vlabel{lem4} Let
$a$ and $b$ tuples in $\calu$ such that $b,\si(a)\in K(a)^{alg}$ and
$\si(b)\in K(b)^{alg}$. 

\begin{enumerate}
\item If $b\in K(a)_{\si^{\pm 1}}$, then
$ld(b/K)\leq ld(a/K)$ and $ild(b/K)\leq ild(a/K)$. 
 \item
  {If $a\in K^{alg}$, then $ld(a/K)=ild(a/K)$.}
\item{If $ld(b/K)=1$, then $K(a,b)_{\si^{\pm 1}}$ is a finite extension of
 $K(a)_{\si^{\pm 1}}$.}
 \item
  {If  some  analysis of $tp(a/K)$ only involves types
non-orthogonal to $\fix(\si)$, then $ld(a/K)=ild(a/K)$.}
\end{enumerate}\end{lem}
\prf (1) Follows from $ld(a,b/K)=ld(a/K)=ld(a/K(b)_{\si^{\pm 1}}) ld(b/K)$, and
similarly for $ild$.

 (2) 
Replacing $a$ by $(a,\ldots,\si^m(a))$ for some $m$, we may assume
    that $ld(a/K)=[K(a,\si(a)):K(a)]$ so that also
    $[K(a,\si(a)):K(\si(a))]=ild(a/K)$. Then
    $[K(a):K]=[K(\si(a)):K]$ gives the result. 

(3) By (2), $ild(b/K(a)_{\si^{\pm 1}})=ld(b/K(a)_{\si^{\pm 1}})=1$; hence, for some $m$,
$K(a,b)_{\si^{\pm 1}}=K(a)_\si(b,\si(b),\ldots,\si^m(b))$.

(4) We know that the limit and inverse limit degrees are multiplicative
    in tower. By (2), $ld(b/K)/ild(b/K)=ld(b/K^{alg})/ild(b/K^{alg})$ is an invariant of the
    extension $K(b)_\si^{alg}/K$. Thus the result holds if $tp(a/K)$ is
    almost internal to $\fix(\si)$: by assumption there are $L=\acl(L)$
    independent from $a$ over $K$ and a tuple $b\in \fix(\si)$ such that
    $\acl(La)=\acl(Lb)$. 
Since $ld(b/L)=ild(b/L)=1$, we have
    $ld(a/L)=ild(a/L)$; because $a$ and $L$ are independent over $K$, we
    have $ld(a/K^{alg})=ld(a/L)$ and $ild(a/K^{alg})=ild(a/L)$. 
By induction on the length of a
    semi-minimal analysis of $tp(a/K)$, we get the result. \qed

\section{From types to isomorphism types}

As defined in \cite{CH}, \cite{CHP}, modularity and fixed-field internality were properties of
a complete type  $tp(c/K)$ in a saturated model   $\calu$ of ACFA.  We
show here that they actually depend only on the difference field
extension $K(c)_{\si^{\pm 1}}/K$.  In the section 2 of the first part
\cite{dyn1}, we present 
the same material differently, defining the notions directly in a way that does not use the embedding.  

We fix two sufficiently saturated  models $\calu$ and
$\calu'$ of ACFA, and
a difference subfield $K$ of $\calu$. Unless otherwise specified, $\tau$ will always be an automorphism
of the form $\si^m\frob^n$, with $m=0$ and $n=1$, or $(m,n)=1$, so that
$SU(\fix(\tau))=1$.  We will usually be
working in $\calu$, when working in $\calu'$ we will indicate it by a
subscript $\calu'$.

\para\vlabel{red}{\bf Reducts}. If $k$ is a positive integer, we denote by $\calu[k]$
 the reduct 
$(\calu,\si^k)$ of the difference field $\calu$. It is also an
existentially closed difference field (Corollary 1.12(1) in \cite{CH}). If $a$ is a tuple in
$\calu$, then we denote by $qftp(a/K)[k]$, $tp(a/K)[k]$, $SU(a/K)[k]$
respectively the
quantifier-free type, type and SU-rank of
the tuple $a$ over $K$ in $\calu[k]$. We will denote by
$\acl_{\si^k}(A)$ the algebraic closure in the sense of $\calu[k]$.

\para\vlabel{lem45}{\bf Codes}. Recall that if $a_1,\ldots,a_m$ are
$n$-tuples in some field, then the code of the
set $\{a_1,\ldots,a_m\}$ is defined  as the tuple $b$
of coefficients of the polynomial 
$\prod_{i=1}^m(X_0+a_i\cdot X)$, where $X=(X_1,\ldots,X_n)$ and
$\cdot$ is the usual dot product of vectors in $n$-space. 

If the tuple $a$ is separably algebraic over the field  $E$, and $a=a_1,a_2,\ldots,a_m$
are the distinct field conjugates of $a$ over $E$, then the code $b$ of
$\{a_1,\ldots,a_m\}$ belongs to $E$. If $a$ is only algebraic over $E$,
then some $p$-th power of $b$ will belong to $E$.

\para \vlabel{remdefi}{\bf The field of definition of the difference
  locus}. Let $L$ be a difference overfield of $K$, and $a$ an $n$-tuple in
$\calu$. We define the {\em difference locus of $a$ over $L$} to be the
smallest  subset of $\calu^n$ defined by difference equations over $L$ and
containing $a$. 

We define the {\em field of definition of the difference
  locus of $a$ over $L$} to be 
 the smallest
difference subfield $E$ of $L\perf$  such that
$E(a)_{\si^{\pm 1}}$ and $L\perf$ are linearly disjoint over $E$, and
denote it by $\qfcb(a/L)$. It can also be described  as the field of definition
of the algebraic locus of the infinite tuple $\{\si^i(a)\mid
i\in\zee\}$ over $L\perf$. 

We 
define the {\em field
of definition of the difference locus of the tuple $a$ over $L/K$} to be 
 the smallest
difference subfield $E$ of $L\perf$ which contains $K$ and is such that
$E(a)_{\si^{\pm 1}}$ and $L\perf$ are linearly disjoint over $E$ and 
denote it by 
$\qfcb_K(a/L)$; note that $\qfcb_K(a/L)=K\qfcb(a/L)$.   

\para\vlabel{rem2defi}{\bf Warning}. Let $a$, $K\subset L$ as above. 
An important observation is that  $\qfcb(a/L)$ is not
  necessarily contained in $L$, but may be purely inseparable over
  $L$. In  {positive} characteristic, it is therefore {\em not to be confused} with the field of definition of the
  $\si$-ideal of difference polynomials over $L$ vanishing at $a$, which
  is contained in $L$, and can be defined as the smallest subfield $E'$ of
  $L$ such that $E'(a)_{\si^{\pm 1}}$ and $L$ are linearly disjoint over
  $E'$. In positive characteristic, the
  two fields of definition are different, it may even happen that $E'$
  is not algebraic over $\qfcb(a/L)$!  

Here is a purely algebraic example where this phenomenon occurs: let $a,b,t$ be algebraically
independent over $\ffi_p$. Let $K=\ffi_p(a^p,b^p,t, a+bt)$, and consider
the point $(a,b)$. Then  $\qfcb(a,b/K)=\ffi_p(a,b)$; however, $K$ is the
field of definition of the ideal $I(a,b/K)=(X^p-a^p,Y^p-b^p, X+tY-(a+bt))$.



\para\vlabel{rem3defi}{\bf More properties of $\qfcb$}. Let $a$,
$K\subset L$ as above, and let $k$ be the prime field. 
\begin{enumerate}

\item 
As in the algebraic case, the field $\qfcb(a/L)$  is contained in the
difference field generated
by realisations of $qftp_{\rm ACF}((\si^i(a)_{i\in\zee})/L)$, where
$qftp_{\rm ACF}$ denotes the type in the reduct to the language of
rings. Also, if $L(a)_{\si^{\pm 1}}\cap L^{sep}=L$ (i.e., if $L(a)_{\si^{\pm 1}}$ is a {\em
  primary} extension of $L$), then there is a unique quantifier-free
type over $L^{alg}$ extending $qftp(a/L)$: this is because
$qftp_{\rm ACF}(\si^i(a)_{i\in\zee})/L)$ is stationary. In that case,
$\qfcb(a/L)$ is contained in the difference field generated by finitely
many realisations of $tp(a/L)$.

\item 
Assume now that $L(a)_{\si^{\pm 1}}/L$ is not primary. Then for some
$\ell$, $\qfcb(a/L)$ is contained in the difference field
generated by finitely many realisations of $\qfcb(a/L)[\ell]$.

\item If $L(a)/L$ is separable, then $\qfcb(a/L)\subset L$.

\item By 5.23.XVIII of \cite{[Co]}, 
$\qfcb(a/L)$ is  finitely  generated  as a ($\si^\ell$-)difference
field, since it is contained in a finitely generated one.
\end{enumerate}

\prf (1) and (2). We will use the following classical facts for
algebraic sets: let $V$ be an absolutely irreducible variety. Then for
some $m$, if $a_1,\ldots,a_m$ are independent generics of $V$,  the
field generated by $a_1,\ldots,a_m$ contains  the
field of definition of $V$. Assume now that $V$ is $L$-irreducible, but
not necessarily absolutely irreducible; 
let $V_0$ be an absolutely irreducible component of $V$, and $k(b)$ its
field of definition. Then the field generated by the field conjugates of
$b$ over $L$ contains the field of definition of $V$, since it contains
the fields of definitions of all irreducible components of $V$.

\smallskip
By 3.8.V of \cite{[Co]}, the topology on cartesian powers of $\calu$ with closed sets the zero-sets  of
difference polynomials over $\calu$ is Noetherian; hence there is an integer $m$ such
that if $(b,\si(b),\ldots,\si^m(b))$ has the same algebraic locus over
$L$ as $(a,\si(a),\ldots,\si^m(a))$, then $a$ and $b$ have the same
difference locus over $L$. This implies that if the tuple $c$ generates
the field 
of definition of the algebraic locus of $(a,\ldots,\si^m(a))$ over $L$,
then $k(c)_{\sipm}=\qfcb(a/L)$. If $L(a)_{\sipm}/L$ is primary, then the
difference locus of $a$ over $L$ is absolutely irreducible, hence,
taking finitely many generic independent realisations of $tp(a/L)$, we
obtain (1).

\smallskip
Let us now prove (2). Let $m$ be as above, and $c$  such that $k(c)$ is
the field of definition of the algebraic locus of
$(a,\si(a),\ldots,\si^m(a))$ over $L^{alg}$. Then the field of
definition of the algebraic locus of $(a,\si(a),\ldots,\si^m(a))$ is
contained in the field generated by the set of (field-) conjugates of
$c$ over $L$. Hence, by the definition of $m$, $\qfcb(a/L)$ is contained in 
the difference field generated by these field conjugates of $c$ over $L$.
By  (1.12) of \cite{CHP}, for
some $\ell$, if $c'$ is a field conjugate of $c$ over $L$, then
$qftp(c'/L)[\ell]=qftp(c/L)[\ell]$. This gives the result.

(3) follows from the fact that $L(a)_{\sipm}$ and $L\perf$ are linearly
disjoint over $L$, and (4) is clear. 
\qed

\bigskip
The following algebraic lemma will be useful

\para\vlabel{lem5}{\bf Lemma}. Let $K\subset L_0\subset L$, $K\subset  M$
and $N$ be fields, such that  $L$
and $M$ are linearly disjoint over $K$,  $L/L_0$ is algebraic, $M/K$ is
regular,  and
$L_0M\subset N\subset LM$. Then $N=L_1M$ where $L_1=L\cap N$
in each of the following cases:
\begin{itemize}
\item [(a)] if $N/L_0M$ is separable, or 
\item [(b)] (char. $p>0$) if $[L_0:L_0^p]= p$ (for instance, if
  $tr.deg(L/K)=1$ and $K$ is perfect).
\end{itemize}

\prf We may assume that $L/L_0$ is finite, since if the conclusion is
false it will be witnessed by a finite subextension.

(a) This is well-known: without loss of generality, $L$ is separable
over $L_0$; let $\hat L$ be the Galois closure of $L$ over $L_0$. Then
the regularity of $M/K$ and the independence of $L$ and $M$ over $K$
imply that $\hat L$ and $M$ are linearly disjoint over $K$; hence 
fields between $L_0M$ and $LM$ correspond to groups between $\gal(\hat
LM/LM)\simeq \gal(\hat L/L)$ and $\gal(\hat LM/L_0M)\simeq \gal(\hat
L/L_0)$. The result follows.

(b)  Using a tower of extensions and the first case, we may assume that
$N/L_0M$ is purely inseparable of degree $p$, and $L/L_0$ is purely
inseparable. Since $[L_0^p:L_0]=p$, 
we get $L\supset L_0^{1/p}$ and  our linear disjointness
assumption implies that $L_0M/L_0$ is separable (as $L$ and $L_0M$ are
linearly disjoint over $L_0$). This implies that 
$N\subseteq L_0^{1/p}M$: assume not, write $N=L_0M(c)$ where $c^p\in L_0M$, and
let $a\in L_0\setminus L_0^p$; if $c\notin L_0^{1/p}M$, then $c,a$ are
$p$-independent, and this implies that $c\notin L_0(a^{1/p^n})$ for any
$n$; i.e., $c\notin L_0\perf M$, which gives us the desired
contradiction. 
\qed

 \begin{lem}  \vlabel{lem1ob0} Let $k\geq 1$, and $(L,\si^k)$ be a
   finitely generated $\si^k$-difference field
  extending $(K,\si^k)$, and such that $L$ is a primary extension of
  $K$. Then there is a $K$-embedding of the $\si^k$-difference field $L$
  into $\calu[k]$.
\end{lem}
 
\prf Since
$\calu[k]$ is also a saturated model of ACFA, it suffices to show the
result for $k=1$. 


The automorphism $\si$ extends uniquely to $K\perf$, and we may
therefore assume that $K$ is perfect. Our primarity hypothesis now implies that
$L$ and $K^{alg}$ are linearly disjoint over $K$, and therefore that
$L\otimes _KK^{alg}$ is a field. Defining $\si$ on $L\otimes_{K}K^{alg}$
so that it extends $\si\rest_{K^{alg}}$ and $\si\rest_L$, we use
\ref{acfa0}(8) to conclude. \qed

\begin{lem}\vlabel{lem2ob0}Let $a$ be a tuple in $\calu$, and assume that $tp(a/K)$ is
  qf-internal to $\fix(\tau)$. Then there is a tuple $c$ of realisations
  of $tp(a/K^{alg})$ such that if $L=K(c)_{\si^{\pm 1}}$, then
  $L(a)_{\si^{\pm 1}}$ is a primary extension of $K(a)_{\si^{\pm 1}}$, $L\dnfo_Ka$, and
  $L(a)_{\si^{\pm 1}}=L(b)_{\si^{\pm 1}}$ for some tuple $b$ in $\fix(\tau)$.
\end{lem}

\prf 
By assumption there is some $L'=\acl(L')$ independent from $a$ over
$K$, and such that $L'(a)_{\si^{\pm 1}}=L'(b)_{\si^{\pm 1}}$ for some tuple $b$ in
$\fix(\tau)$. 

Let $(a_i,b_i)_{i\in\nat}$ be an independent sequence of realisations of
$tp(a,b/L')$, with $(a_0,b_0)=(a,b)$. Then for some $n$,
$K(a_1,b_1,\ldots,a_n,b_n)_{\si^{\pm 1}}$ contains $\qfcb_K(a,b/L')$. Let
$L=K(a_1,\ldots,a_n)_{\si^{\pm 1}}$. Then
$L(a,b_1,\ldots,b_n)_{\si^{\pm 1}}=L(b,b_1,\ldots,b_n)_{\si^{\pm 1}}$, i.e.,
$L(a)_{\si^{\pm 1}}\subset L\fix(\tau)$. Observe that $L$ and $K(a)_{\si^{\pm 1}}^{sep}$ are
linearly disjoint over $K(a)_{\si^{\pm 1}}\cap K^{sep}$: this is because the $a_i$'s
are independent realisations of $tp(a/L')$, and in particular of
$tp(a/K^{alg})$, and because $a\dnfo_KL'$. Hence $L(a)_{\si^{\pm 1}}$ is a 
primary extension of $K(a)_{\si^{\pm 1}}$, and is generated over $K$ by
realisations of $tp(a/K^{alg})$.   \qed

\begin{prop}\vlabel{ob0}
Let $a\in\calu$ and $\ell$ a positive integer, $\fix(\tau)$ a fixed
field. 
\begin{enumerate}
\item Assume that $tp(a/K)$ is
qf-internal to  $\fix(\tau)$. Also, assume that
$\varphi:K(a)_{\si^{\pm 1}}\to K'(a')_{\si^{\pm 1}}\subset\calu'$ is an isomorphism with $\varphi(K)=K'$, $\varphi(a)=a'$. Then $tp_{\calu'}(a'/K')$ is
qf-internal to $\fix(\tau)$. 

\item $tp(a/K)[\ell]$ is qf-internal to
  $\fix(\tau^\ell)$ if and only if $tp(a/K)$ is qf-internal to
  $\fix(\tau)$. 
\end{enumerate}

\end{prop}

\prf (1) Apply Lemmas \ref{lem2ob0} and \ref{lem1ob0} (with $k=1$).

(2) Assume first that $tp(a/K)$ is qf-internal to $\fix(\tau)$. Then
    $tp(a/K)[\ell]$ is qf-internal to $\fix(\tau^\ell)$ because
    $\fix(\tau)\subset \fix(\tau^\ell)$: by hypothesis there is
    $L=\acl(L)$ which is independent from $a$ over $K$ and such that
    $a\in L\fix(\tau)\subset L\fix(\tau^\ell)$.  

Assume now that $tp(a/K)[\ell]$ is qf-internal to $\fix(\tau^\ell)$. Then so
are $tp(\si^i(a)/K)[\ell]$ for all $i$, and, replacing $a$ by
$(a,\si(a),\ldots,\si^{m}(a))$ for some $m$, we may therefore assume that
$\si(a)\in K(a)^{alg}$, and $K(a)_{\si^{\pm 1}}=K(a)_{\si^{\pm\ell}}$. 
By \ref{lem2ob0}, in
$\calu[\ell]$, there is a $\si^\ell$-difference field $L$ such that
$L(a)_{\si^{\pm \ell}}\subset L\fix(\tau^\ell)$, 
$L\dnfo_K a$, and $L(a)_{\si^{\pm\ell}}/K(a)_{\si^{\pm \ell}}$ is primary. 
 However,
it may be that the $\si$-difference subfield of $\calu$ generated by $L$
is not independent from $a$ over $K$. Since
$\si$ and $\si^\ell$ extend uniquely to the perfect closure of a field, we may assume
that $L$ and $K$ are perfect, so that now $L(a)_{\si^{\pm\ell}}$ is a regular
extension of $K(a)_{\si^{\pm\ell}}=K(a)_{\si^{\pm 1}}$. For each $i=0,\ldots,\ell-1$, we
choose $L_i$ realising $\si^i(tp(L/K(a)_{\si^{\pm 1}})[\ell])$ and such that $L_0=L$,
$L_1,\ldots,L_{\ell-1}$ are independent over $K(a)_{\si^{\pm\ell}}=K(a)_{\si^{\pm 1}}$. Since they
are isomorphic copies of the regular extension
$L(a)_{\si^{\pm\ell}}/K(a)_{\si^{\pm 1}}$, they are in fact linearly disjoint over
$K(a)_{\si^{\pm 1}}$. Hence, reasoning as in 1.12 of \cite{CH}, we can define an
extension $\rho$ of $\si$ on the composite field $L_0\cdots L_{\ell-1}$,
such that $\rho^\ell$ coincides with $\si^\ell$ on $L_0$. Thus, by Lemma
\ref{lem1ob0}, there is a $K(a)_{\si^{\pm 1}}$-embedding $\varphi$ of the $\si$-difference
field $(L_0\cdots L_{\ell-1},\rho)$ into $\calu$; then $a\in
\varphi(L)\fix(\tau^\ell)$, and we are done. \qed

\begin{lem}\vlabel{lem3ob0}Let $a$ be a tuple in $\calu$, and assume
  that $tp(a/K)$ is non-orthogonal to $\fix(\tau)$. Then there is $e\in
  K(a)_{\si^{\pm 1}}$ such that $tp(e/K)$ is qf-internal to $\fix(\tau)$.\end{lem} 

\prf By assumption (and Fact \ref{lem2}) there are tuples $b$ and $c$, with $c\dnfo_Ka$, $b\in
\fix(\tau)$, and $b\in K(c,a)_{\si^{\pm 1}}$, $b\notin K(c)_{\si^{\pm 1}}^{alg}$. Reason
exactly as in the proof of \ref{lem2ob0} to show that $c$ can 
be chosen so that $K(c,a)_{\si^{\pm 1}}$ is a primary extension of
$K(a)_{\si^{\pm 1}}$. As $b\in K(a,c)_{\sipm}$,
it follows that $\qfcb_K(b,c/K(a)_{\si^{\pm 1}})$ is contained in the
difference field generated by independent realisations of
$tp(b,c/K(a)_{\si^{\pm 1}})$, i.e., if $K(e)_{\si^{\pm 1}}=\qfcb_K(b,c/K(a)_{\si^{\pm 1}})$, then
$tp(e/K)$ is qf-internal to $\fix(\tau)$. For some $n$, $e^{p^n}\in
K(a)_{\si^{\pm 1}}$. \qed

\begin{prop}\vlabel{ob0b} Let $a$ be a tuple in $\calu$, and $\ell$ a
  positive integer. 
\begin{enumerate}
\item $tp(a/K)$ is one-based if and only if $tp(a/K)[\ell]$ is
  one-based.
\item If $\varphi:K(a)_{\si^{\pm 1}}\to K'(a')_{\si^{\pm 1}}\subset \calu'$ is an
  isomorphism  with $\varphi(K)=K'$, $\varphi(a)=a'$, and
if  $tp(a/K)$ is one-based, then so is $tp_{\calu'}(a'/K')$. 
\item If $tp(a/K)$ is non-orthogonal to a one-based type, then there is
  $b\in K(a)_{\si^{\pm 1}}$ such that $tp(b/K)$ is one-based.
\end{enumerate}
\end{prop}

\prf (1)
Assume first that $tp(a/K)[\ell]$ is one-based,  let $\bar a$ be a
    tuple of realisations of $tp(a/K)$, and $L=\acl(L)\supset K$. Observe
    that for each $i$, $tp(\si^i(a)/K)[\ell]$ is also one-based, and
    therefore so is $tp(\bar a,\si(\bar a),\ldots,\si^{\ell-1}(\bar a)/K)[\ell]$. Hence 
    $\acl_{\si^{\pm\ell}}(K,\bar a,\si(\bar a),\ldots,\si^{\ell-1}(\bar
    a))=\acl(K\bar a)$ and $L$ are linearly disjoint over their
    intersection. This 
     shows that $tp(a/K)$ is one-based.

Assume that $tp(a/K)$ is one-based, but not $tp(a/K)[\ell]$. Replacing
$a$ by $(a,\si(a),\ldots,\si^m(a))$ for some $m$, we may assume that
$\si(a)\in K(a)^{alg}$ and $K(a)_{\si^{\pm 1}}=K(a)_{\si^{\pm\ell}}$.  Then, using the
semi-minimal analysis in $\calu[\ell]$ and the trichotomy, there
is $c\in K(a)^{alg}$ such that $tp(c/K)[\ell]$ is one-based, and
$tp(a/K(c)_{\si^{\pm\ell}})[\ell]$ is non-orthogonal to $\fix(\tau^\ell)$ for
some fixed field $\fix(\tau)$. Then also each
$tp(\si^i(c)/K)[\ell]$ is one-based, and therefore
$tp(a/K(c)_{\si^{\pm 1}})[\ell]$ is non-orthogonal to $\fix(\tau^\ell)$ as
well (see \ref{onebased}). Enlarging $c$, we may assume that $K(c)_{\si^{\pm 1}}=K(c)_{\si^{\pm\ell}}$. By Lemma \ref{lem3ob0}, there is a tuple $b\in K(c,a)_{\si^{\pm\ell}}$ such
that $tp(b/K(c)_{\si^{\pm\ell}})[\ell]$ is qf-internal to
$\fix(\tau^\ell)$. Since $K(c)_{\si^{\pm 1}}=K(c)_{\si^{\pm\ell}}$, Lemma~\ref{lem2ob0}
gives $tp(b/K(c)_{\si^{\pm 1}})$ qf-internal to $\fix(\tau)$, which contradicts
the one-basedness of $tp(a/K)$.

(2) We may assume that $\si(a)\in K(a)^{alg}$. The proof is by induction
on $tr.deg(K(a)_{\si^{\pm 1}}/K)$. If it is $0$, there is nothing to prove. Assume
that $tp_{\calu'}(a'/K')$ is not one-based. 
Then there is $c'\in
K'(a')^{alg}$ such that $tp_{\calu'}(c'/K')$ is one-based, and $tp_{\calu'}(a'/K(c')_{\si^{\pm 1}})$ is
non-orthogonal to some fixed field $\fix(\tau)$. We may assume that
$\si(c')\in K'(c')^{alg}$, and that either $c'\notin {K'}^{alg}$, or $c'\in
K'$.

If $c'\in K'$, then Lemma \ref{lem3ob0} and Proposition
\ref{ob0} imply $tp(a/K)$ non-orthogonal to $\fix(\tau)$, a
contradiction.

Assume that  $c'\notin K'$. 
By 1.12(3) of
\cite{CHP}, there are  $\ell\geq 1$ and $c\in\calu$ such that $\varphi$ extends to a
$\si^\ell$-isomorphism $\varphi:K(a)_{\si^{\pm 1}}(c)_{\si^{\pm\ell}}\to
K(a')_{\si^{\pm 1}}(c')_{\si^{\pm\ell}}$.
Since $\si(c')\in K'(c')^{alg}$, 
$tp_{\calu'}(a'/K(c')_{\si^{\pm\ell}})[\ell]$ is   non-orthogonal to
 $\fix(\tau^\ell)$. By Lemma \ref{lem3ob0} and Proposition \ref{ob0}, so is
 $tp(a/K(c)_{\si^{\pm\ell}})[\ell]$. As
$tr.deg(K'(c')_{\si^{\pm 1}}/K')<tr.deg(K'(a')_{\si^{\pm 1}}/K')$, by induction
hypothesis, we know that $tp(c/K)[\ell]$ is one-based, as are all
$tp(\si^i(c)/K)[\ell]$. 
Hence, $tp(a/K(c)_{\si^{\pm 1}})[\ell]$ is non-orthogonal
to $\fix(\tau^\ell)$. We then reason as in the previous case to get a
contradiction.

(3) Our assumption implies that there is $c\in \acl(Ka)$ with $tp(c/K)$
one-based (see \ref{ana}). Let $c=c_1,c_2,\ldots, c_m$ be the
distinct field conjugates of $c$ over $K(a)_{\si^{\pm 1}}$. By 1.12 of
\cite{CHP}, there is $\ell\geq 1$ such that they all have the same
quantifier-free type over $K$ in $\calu[\ell]$. By (1), $tp(c/K)[\ell]$ is
one-based, and by (2), so are $tp(c_i/K)[\ell]$ for $i\geq 1$. Hence if
$b\in K(a)_{\si^{\pm 1}}$ is (some $p^n$-power of) the code of the set
$\{c_1,\ldots,c_m\}$, then $tp(b/K)[\ell]$ is one-based, and so is
$tp(b/K)$. \qed

\para{\bf Definition}. Let $a$ be a tuple in $\calu$. We say that
$K(a)_{\si^{\pm 1}}/K$ is {\em primitive} if $tr.deg(K(a)_{\si^{\pm 1}}/K)<\infty$, and whenever $b\in K(a)_{\si^{\pm 1}}$, then either
$b\in K^{alg}$ or $a\in K(b)_{\si^{\pm 1}}^{alg}$.  (This agrees with the
definition of primitive algebraic dynamics given in \cite{dyn1}.)

The previous results \ref{lem3ob0}, \ref{ob0} and \ref{ob0b} then   immediately yield:

\begin{prop}\vlabel{lemcor2}Assume that $K(a)_{\si^{\pm 1}}/K$ is
  primitive and of finite transcendence degree (over $K$). Then
\begin{enumerate}
\item either $tp(a/K)$ is one-based,
\item or there is $d\in K(a)_{\si^{\pm 1}}$ such that $tp(d/K)$ is qf-internal
  to some fixed field $\fix(\tau)$ and $a\in K(d)_{\si^{\pm 1}}^{alg}$.
\item Whether (1) or (2) holds only depends on the isomorphism type of
  the difference field extension $K(a)_{\si^{\pm 1}}/K$, i.e.: let $\varphi: K(a)_{\si^{\pm 1}}\to K'(a')_{\si^{\pm 1}}\subset \calu'$ be an isomorphism with
  $\varphi(K)=K'$, $\varphi(a)=a'$. 
Then $tp(a/K)$ is one-based [resp. $tp(d/K)$ is qf-internal to
  $\fix(\tau)$ and $a\in K(d)_{\si^{\pm 1}}^{alg}$] if and only if
  $tp_{\calu'}(a'/K')$ is one-based [resp. if $d'=\varphi(d)$, then
  $tp_{\calu'}(d'/K')$ is qf-internal to $\fix(\tau)$ and $a'\in K'(d')_{\si^{\pm 1}}^{alg}$]. 
\end{enumerate}
\end{prop}

\bigskip

Clearly, given some tuple $a$, one can find tuples $a_1,\ldots,a_n\in
K(a)_{\si^{\pm 1}}$ such that $K(a)_{\si^{\pm 1}}=K(a_1,\ldots,a_n)_{\si^{\pm 1}}$, and for each $i$,
$K(a_1,\ldots,a_i)_{\si^{\pm 1}}/K(a_1,\ldots,a_{i-1})_{\si^{\pm 1}}$ is either 
primitive  or algebraic.
Thus, Proposition \ref{lemcor2} has the following immediate consequence:

\begin{thm}\vlabel{cor2}\vlabel{propo4}  Assume that $a$ has
finite SU-rank 
over $K$. Then there are $a_1,\ldots,a_n$ such that
$K(a_1,\ldots,a_n)_{\si^{\pm 1}}=K(a)_{\si^{\pm 1}}$, and 
 for each $i$,
$tp(a_i/K(a_1,\ldots,a_{i-1})_{\si^{\pm 1}})$ is of one of the following three
kinds:
\begin{itemize}
\item[\rm (a)] algebraic,
\item[\rm (b)] one-based, 
\item[\rm (c)] qf-internal to $\fix(\tau)$ for some $\tau$.
\end{itemize}
\end{thm}

\para These  results have many easy consequences. Let us mention two,
others can be derived in a similar way. We fix a subset $\pi$ of
$$\{(\si^n(x^{p^m})=x: (n,m)=1, n>0\}\cup \{ \{\hbox {all one-based
  types of SU-rank }1\}\}.$$

\begin{prop}\vlabel{propo5}Let $a$ be a tuple, $\ell\geq 1$, $E$ a
  difference subfield of $\calu$, and assume
  that $tp(a/K)$ 
  is $\pi$-analysable. 
\begin{enumerate}
\item Then so is $qftp(a/K)[\ell]$.
\item If $a\in E^{alg}$, and $b$ is a code for the set
  of field conjugates of $a $ over $E$, then so is  $tp(b/K)$.
\item Let $\varphi:K(a)_{\si^{\pm 1}}\to K'(a')_{\si^{\pm 1}}\subset \calu'$ be an
  isomorphism with $\varphi(K)=K'$, $\varphi(a)=a'$. Then
  $tp_{\calu'}(a'/K')$ is $\pi$-analysable. 

\item  Statements (1), (2) and (3) hold if one replaces $\pi$-analyzable by
  qf-internal to $\pi$. 

\end{enumerate}
 \end{prop}

\prf (1) Let $(a_1,\ldots,a_n)\in K(a)_{\si^{\pm 1}}$ satisfy the conclusion of Theorem
\ref{cor2} for the extension $K(a)_{\si^{\pm 1}}/K$. 
Then by Propositions \ref{ob0} and \ref{ob0b} so does the tuple
$(a_1,\si(a_1),\ldots,\si^{\ell-1}(a_1),a_2,\si(a_2),\ldots,\si^{\ell-1}(a_n))$
for the extension $K(a,\si(a),\ldots,\si^{\ell-1}(a))_{\si^{\pm\ell}}/K$ in
$\calu[\ell]$. As $K(a)_{\si^{\pm\ell}}\subset K(a)_{\si^{\pm 1}}$, this gives the
result. 

(2) By 1.12 of \cite{CHP}, for some $\ell$, all field conjugates of
$a$ over $E$ satisfy the same quantifier-free type in $\calu[\ell]$. The
result follows from (1).

(3) Clear from Theorem \ref{cor2} and Proposition \ref{lemcor2}, 

(4) Identical proof. \qed


\begin{thm} \vlabel{theo3} 
 Let $K\subset L$ be difference fields, and $a$ a
tuple. Let $\pi$ be a set of types as above. Let $b\in L\perf$ be such that
$K(b)_{\si^{\pm 1}}=\qfcb_K(a/L)$. Assume that there is some $c\in K(a)_{\si^{\pm 1}}$ which is
independent from $L$ 
over $K$ and such that $tp(a/K(c)_{\si^{\pm 1}})$ is $\pi$-analysable. Then so is
$tp(b/K)$. \end{thm}

\prf 
Let $d$ be such that 
$K(d)_{\si^{\pm 1}}=\qfcb_K(a/L^{alg})$. As $c\in K(a)_{\si^{\pm 1}}$, $d$ is
contained in the difference field generated over $K$ by a tuple
$(c^1,a^1)$ of finitely many
$L$-independent realisations of $tp(c,a/L^{alg})$. As $c\dnfo_KL$, it
follows that 
$tp(d/K)$ is qf-internal to the set of $\Aut(\calu/K)$-conjugates of $tp(a/K(c)_{\si^{\pm 1}})$,
and is therefore $\pi$-analysable. 

If $d'$ is a field conjugate of $d$ over $L$, then $d'$ realises
$qftp(d/L)[\ell]$ for some $\ell\geq 1$.  
As $b$ belongs to the field generated
over $K$ by field conjugates of $d$ over $L$, i.e., by realisations of
$qftp(d/L)[\ell]$ for some $\ell$ (see the remark at the end of
\ref{remdefi}), Proposition
\ref{propo5} implies that $tp(b/K)$ is $\pi$-analysable. \qed

\section{Descent}

\begin{prop}\vlabel{theo1}
 Let $(\calu,\si)$ be a  
  difference field, and $K_1\subset K_2$ be subfields of the fixed field
  $\fix(\si)$ of $\calu$,
  with $K_2/K_1$ regular and $SU(a/K_2)<\infty$.  
   Let
  $a,b$ be tuples in $\calu$ 
  such that: \begin{itemize}
\item[\rm (a)] $a\in K_2(b)_{\si^{\pm 1}}$;
\item [\rm (b)]$K_1(b)_{\si^{\pm 1}}$ is linearly disjoint from
  $K_2$ over $K_1$; 
\item [\rm (c)]$tp(a/K_2)$ is hereditarily orthogonal to
  $\fix(\si)$. 
\end{itemize} 

Then there is a tuple $c $ in $K_1(b)^{\rm perf}_{\si^{\pm 1}}\cap
K_2(a)\perf_{\si^{\pm 1}}$ such that $a\in K_2(c)_{\si^{\pm 1}}$. If $tr.deg(K_2(a)_\si/K_2)=1$ and
$K_1$ is perfect, or if the characteristic is $0$, 
then $K_1(b)_{\si^{\pm 1}}\cap K_2(a)_{\si^{\pm 1}}=K_2(c)_{\si^{\pm 1}}$.

\end{prop}
\prf  We may assume $(\calu,\si)$ is existentially closed and saturated. 
Let $e$ be a finite tuple in $K_2$ be such that $a\in K_1(e,b)_{\si^{\pm 1}}$, and
$K_1(e,a)_{\si^{\pm 1}}$ and $K_2$ are linearly disjoint over
$K_1(e)_{\si^{\pm 1}}$. 
Let 
$K_1(d)_{\si^{\pm 1}}=\qfcb_{K_1}(e,a/K_1(b)_{\si^{\pm 1}})$. 

As 
$K_2/K_1$ is separable (and $b\dnfo_{K_1}K_2$, $e\in K_2$), so is the extension 
$K_1(b,e,a)_{\sipm}=K_1(b,e)_{\si^{\pm 1}}$ of $K_1(b)_{\si^{\pm 1}}$; hence $d\in
K_1(b)_{\si^{\pm 1}}$. From the linear disjointness of $K_1(b)_{\sipm}$
and $K_1(a,d,e)_{\sipm}$ over $K_1(d)_{\sipm}$ and the fact that $a\in
K_1(b,e)_{\sipm}$, we obtain $a\in K_1(d,e)_{\sipm}$. 
Since
$e\dnfo_{K_1}b$, 
by Theorem 
\ref{theo3}, $tp(d/K_1)$ is hereditarily orthogonal to $\fix(\si)$.

Let $K_1(c)_{\si^{\pm 1}}=\qfcb_{K_1}(d/K_1(e,a)_{\si^{\pm 1}})$; 
by Theorem \ref{theo3}, $tp(c/K_1)$ is hereditarily
orthogonal to 
$\fix(\si)$, which implies $c\dnfo_{K_1}e$ (since $e\in K_2\subset \fix(\si)$). By definition of $c$,
the fields
$K_1(c,d)_{\si^{\pm 1}}$ and $K_1(e,a)^{\rm perf}_{\si^{\pm 1}}$ are linearly disjoint over
$K_1(c)_{\si^{\pm 1}}$; as $a\in K_1(e,d)_{\si^{\pm 1}}$, this implies that
$a\in K_1(c,e)_{\si^{\pm 1}}$, and therefore $K_2(a)^{\rm perf}_{\si^{\pm 1}} =
K_2(c)^{\rm perf}_{\si^{\pm 1}}$. 

For the first assertion, it remains to show that $c$ can be taken in
$K_1(b)^{\rm perf}_{\si ^{\pm 1}}$.

We now use that $K_2$ is a regular extension of $K_1$.  We know that $c\in
K_2(a)^{\rm perf}_{\si^{\pm 1}} \subset K_2(b)^{\rm perf}_{\si^{\pm
    1}}$; we also know that $tp(c/K_1)$ is 
hereditarily orthogonal to $\fix(\si)$; as $b\dnfo_{K_1}K_2$ and $c\in
\acl(K_2b)$, this implies $c\in \acl(K_1b)$. The regularity of $K_2/K_1$ then implies
that $K_2\perf$ and 
$K_1(b,c)^{\rm perf}_{\si^{\pm 1}}$ are linearly disjoint over $K_1\perf$, and gives
$c\in K_1(b)^{\rm perf}_{\si^{\pm 1}}$. 

Note that if the characteristic is $0$, then we directly obtain that
$K_1(c)_{\si^{\pm 1}}=K_1(b)_{\si^{\pm 1}}\cap K_2(a)_{\si^{\pm
    1}}$. Assume that the characteristic 
is positive, that $K_1$ is perfect and that $tr.deg(K_2(a)_{\si^{\pm
    1}}/K_2)=1$. Then $tr.deg(K_1(c)_{\si^{\pm 1}}/K)=1$ and we can
apply Lemma \ref{lem5}, as for some $n$, we have $c^{p^n}\in
K_2(a)_{\si^{\pm 1}}$:  if $c'$ is such that $K_1(c)_{\si^{\pm 1}}\cap K_2(a)_{\si^{\pm 1}}
=K_1(c')_{\si^{\pm 1}}$, then $K_2(c')_{\si^{\pm 1}}=K_2(a)_{\si^{\pm
    1}}$.

\begin{thm}\vlabel{theo1a}
 Let $(\calu,\si)$ be a  
  difference field, and $K_1\subset K_2$ be subfields of the fixed field
  $\fix(\si)$ of $\calu$,
  with $K_2/K_1$ regular.  
   Let
  $a,b$ be tuples in $\calu$ 
  such that  $SU(a/K_2)<\infty$ and: \begin{itemize}
\item[\rm (a)] $a\in K_2(b)_{\si}$;
\item [\rm (b)]$K_1(b)_{\si}$ is linearly disjoint from
  $K_2$ over $K_1$; 
\item [\rm (c)]$tp(a/K_2)$ is hereditarily orthogonal to
  $\fix(\si)$. 
\end{itemize} 

Then there is a tuple $c $ in $K_1(b)\perf_\si$ such that $a\in K_2(c)_\si$ and $c$ is
purely inseparable over $K_2(a)_\si$. If $tr.deg(K_2(a)_\si/K_2)=1$ and
$K_1$ is perfect, 
then one can choose $c$ so that  $K_1(b)_\si\cap K_2(a)_\si=K_1(c)_\si$. 

\end{thm}

\prf Let $c$ be given by  Proposition \ref{theo1}; we may choose it such
that $K_1(c)_{\si^{\pm
    1}}\subseteq K_1(c)_\si^{alg}$ and $a\in K_2(c)_\si$ (replace $c$ by
 $(\si^{-m}(c),\ldots,\si^{-n}(c))$ for some $n\leq m$).  Choose 
$i\geq
0$ such that  $\si^i(c)\in K_1(b)^{\rm perf}_\si\cap K_2(a)^{\rm
  perf}_\si$.

Choose a tuple $d$ such that $K_1(d)_{\si}=K_1(c)_\si
\cap 
K_1(\si^i(c))_\si^{sep}$. Then 
$a\in K_2\acl(K_1d)$. Also, we know that $\si^i(a)\in
K_2(d)_\si$, and this implies that $K_2(d,a)_\si$ is a finite extension
of $K_2(d)_\si$, which is contained in $K_2(c)_\si$. As
$K_2/K_1$ is regular, $(K_1(b)_\si)^{alg}$ and $K_2$ are
linearly disjoint over $K_1$; by Lemma \ref{lem5}(a), this gives
$K_2(d,a)_\si\cap K_2(d)^{sep}_\si=K_2(d)_\si$, and therefore $a\in
K_2(d)^{\rm perf}_\si$. Replacing $d$ by $d^{p^{-m}}$ for some $m$, we
have $a\in K_2(d)_\si$. This gives us the general case.

Assume now that $K_1$ is perfect and $tr.deg(K_2(a)_\si/K_2)=1$; then
also $tr.deg(K_1(c)_\si/K_1)=1$. We reason as before: the element $c$
can now be chosen in $K_2(a)_{\si^{\pm 1}}\cap K_1(b)_{\si^{\pm 1}}$,
and such that $a\in K_2(c)_\si$, and for some $i$, $\si^i(c)\in
K_2(a)_\si$. Then 
$K_2(a)_\si$ is a field between $K_2(c)_\si$ and
$K_2(\si^i(c))_\si$, and  Lemma
\ref{lem5} now  gives us a tuple $d$ such that $K_1(d)_\si=K_2(a)_\si\cap
K_1(c)_\si$, and $K_2(a)_\si=K_2(d)_\si$.

\begin{thm}
\vlabel{theo2} Let $K_1\subset K_2$ be fields, with $K_2/K_1$
  regular,  $V_2$ an irreducible variety over $K_2$, and
  $(V_2,\phi_2)\in AD_{K_2}$. Assume that $(V_2,\phi_2)$ is
  primitive and $\deg(\phi_2)>1$.
  Assume furthermore that for some $n\geq 1$,  $(V_2,\phi_2^n)$ is dominated
  (in $AD_{K_2}$)
  by some object of $AD_{K_1}$. 

\begin{enumerate}
\item There is some variety $V_3$ defined over $K_1$, and a
  dominant constructible map $\phi_3:V_3\to V_3$ also defined over $K_1$, a constructible
  isomorphism $h:(V_2,\phi_2)\to (V_3,\phi_3)$. 
\item Assume that the characteristic is $0$, or that $K_1$ is perfect
  and $\dim(V_2)=1$. Then $(V_2,\phi_2)$ is rationally isotrivial, i.e., there is
  some $(V_3,\phi_3)\in
AD_{K_1}$ which is isomorphic to $(V_2,\phi_2)$ (in $AD_{K_2}$).
\end{enumerate}
\end{thm}

\prf  
Let $(V_1,\phi_1)\in
AD_{K_1}$ dominate $(V_2,\phi^n_2)$, via a map $f$. The map $\phi_1$ permutes the absolutely irreducible
components of $V_1$, and for some $m$, $\phi_1^m$ leaves them
invariant. 

Let $\calu$ be a saturated, existentially closed   
difference field containing $K_2$ and such that $\si$ is the identity on
$K_2^{alg}$. By saturation of $\calu$, there is a generic $b$ of $V_1$
over $K_2$ which satisfies $\si(b)=\phi_1^m(b)$. 
 By assumption, $f(b)=a$ is a generic of
$V_2$ and satisfies $\si(a)=\phi_2^{mn}(a)$.

\medskip\noindent
{\bf Claim}.  $tp(a/K_2)$ is hereditarily orthogonal to $\fix(\si)$. 

By primitivity and
Proposition \ref{propo4}, if the claim is false then $tp(a/K_2)$ is
almost internal to $\fix(\si)$. By Lemma \ref{lem4}(4), we have
$ild(a/K_2)=ld(a/K_2)=1$, which contradicts the assumption.

\medskip

 By Theorem \ref{theo1a}, there is $c\in K_1(b)^{\rm perf}\cap
 K_2(a)^{\rm perf}$ such that $a\in K_2(c)$. Thus there is a
 constructible isomorphism $g:V_2\to V_3$, where $V_3$ is the algebraic
 locus of $c$ over $K_1$, and such that $g(a)=c$. 
%
%
%
 Consider now
$a'=\phi_2(a)$, and $c'=g(a')$. Then
$\si(a')=\phi_2^{mn+1}(a)=\phi_2^{mn}(a')$, and, because $V_2$ was
$K_2$-irreducible, we have 
$$qftp(a/K_2)=qftp(a'/K_2), \qquad qftp(c/K_2)=qftp(c'/K_2).$$
As $tp(c/K_2)$ is hereditarily orthogonal to $\fix(\si)$ (i.e., to
$tp(K_2/K_1)$) and $c\dnfo_{K_1}K_2$, it follows (by \ref{theo3}) that
$cc'\dnfo_{K_1}K_2$. As $K_2/K_1$ is regular, the fields $K_1(c,c')$ and
$K_2$ are linearly disjoint over $K_1$, and this implies that $c'\in
K_1(c)\perf$ (recall that $c'\in K_2(c)\perf$). Hence, for some
dominant constructible $\phi_3:V_3\to V_3$ sending $c$ to $c'$, we have
$g:(V_2,\phi_2)\simeq (V_3,\phi_3)$. This shows (2).

Assume now that $K_1$ is perfect, and $\dim(V_2)=1$. Using Theorem
\ref{theo1a}, we get $c$ such that $K_1(b)\cap K_2(a)=K_1(c)$ and
$K_2(c)=K_2(a)$, whence a birational isomorphism $g:V_2\to
V_3$. Defining $a'=\phi_2(a)$ and $c'=g(a')$, we obtain $c'\in
K_2(c)\cap K_1(c)^{alg}=K_1(c)$, i.e., if $\phi_3=g\phi_2g\inv$, then
$\phi_3$ is a dominant morphism defined over $K_1$, and
$g:(V_2,\phi_2)\simeq (V_3,\phi_3)$. \qed

\begin{rem}\begin{enumerate}
\item The proof of Theorem \ref{theo2} actually shows more: in the above notation,
the fact that $c\in K_1(b)\perf_{\si}$ means that the composed
constructible morphism
$g\circ f:V_1\to V_2\to V_3$
is defined over $K_1\perf$. 
\item The hypothesis on the degree of $\phi$ can be replaced by
  $(V,\phi)$ fixed-field-free.
\item Inspection of the proof of Proposition \ref{theo1} shows that the
  hypotheses that 
\begin{quote}\em $K_1\subset K_2\subset \fix(\si)$ and $tp(a/K_2)$ is hereditarily
  orthogonal to $\fix(\si)$ \end{quote}
can be be changed to:
\begin{quote}\em $K_2/K_1$ is regular, and $tp(K_2/K_1)$ is hereditarily orthogonal to all types
  occurring in a semi-minimal analysis of $tp(a/K_2)$.\end{quote}
\item Also, as observed (and proved) in \ref{1.4} of \cite{dyn1}, to obtain
 isogeny isotriviality, the only assumption needed on
$K_1\subset K_2$ and $a$ is that $tp(a/K_2)$ is one-based. 
\item Theorem \ref{theo2}, together with Lemma \ref{limited2} and 
    Remark \ref{limited3}(2h) in \cite{dyn1}, gives the result of M. Baker, Theorem 1.6  in
  \cite{baker}. 

\end{enumerate}
\end{rem}

\def\proj{\Pp}
\para\vlabel{counterex}{\bf Example}.    {\em An algebraic dynamics over
  $k(t)$, dominated by an 
algebraic dynamics over 
$k$, but not isogenous to a difference variety over $k$. }   

\smallskip\noindent
Let $H$ be a vector extension of an Abelian variety; i.e. there
exists an exact sequence of algebraic groups $0 \to V \to H \to A \to 0$ with $A$ an Abelian 
variety, $V \cong \ga^n$ a vector group.  Assume $\dim(V) = 2$ and 
$\Hom(H,\ga) = (0)$.  Assume $H$ is defined over $k$, and identify the projective
space $\proj V $ with  $\proj^1$.   In particular a transcendental element $t$ gives a one-dimensional subspace $V_t$ of $V$.   Let $H_t = H / V_t$.  Fix $h \in H$ generating $H$;
it suffices that the image of $h$ in $A$ generate $A$; if $A$ is a simple Abelian variety, it
suffices therefore that the image is non-torsion.  Let $Y=  (H,T(h))$ and $X_t = (H_t, T(h_t))$
where $h_t$ is the image of $h$ in $H_t$, and $T(g)$ denotes translation by $g$.  Then $Y$ dominates $X_t$.  But $X_t$
is clearly not isotrivial.   


\bigskip \noindent
Current addresses:

\bigskip\noindent
UFR de Math\'ematiques \par\noindent
Universit\'e Paris 7 - Case 7012
\par\noindent 2, place Jussieu
\par\noindent
75251 Paris Cedex 05 \par\noindent
France \par\noindent
e-mail: {\tt zoe@logique.jussieu.fr}

\bigskip\noindent
Institute of Mathematics \par\noindent
The Hebrew University \par\noindent
Givat Ram \par\noindent
91904 Jerusalem \par\noindent
Israel \par\noindent
e-mail: {\tt ehud@math.huji.ac.il}


\bigskip\noindent


\begin{thebibliography}{99}

\itemsep=\smallskipamount
\bibitem{baker}{M. Baker, A finiteness theorem for canonical heights attached
to rational maps over function fields, available at {\tt
  arXiv:math.NT/0601046}.}









\bibitem{CH}{Z. Chatzidakis, E. Hrushovski, Model theory of difference 
fields,
Trans. Amer. Math. Soc. 351 (1999), pp. 2997-3071.}
\bibitem{CHP}{Z. Chatzidakis, E. Hrushovski, Y. Peterzil, Model
theory of difference fields, II: Periodic ideals and the trichotomy
in all characteristics, Proc. London Math. Soc. (3) 85 (2002) 257 --
311.}

\bibitem{dyn1}Z. Chatzidakis, E. Hrushovski, Difference fields and
  descent in 
  algebraic dynamics, I, this volume.
\bibitem{dyn3}Z. Chatzidakis, E. Hrushovski, Difference fields and descent in
  algebraic dynamics, III, in preparation.


\bibitem{[Co]}{R.M. Cohn, {\em Difference algebra}, Tracts in
Mathematics 17, Interscience Pub. 1965.}







 \bibitem{h-emm}
 Ehud Hrushovski,  The Manin-Mumford conjecture and the model theory of difference
 fields, Annals of Pure and Applied Logic, vol. 112 (2001), no. 1,
 pp. 43 -- 115.




    




 
  

\bibitem{[W]}{F. Wagner, Some remarks on one-basedness, J. of Symb. Logic 69
Nr 1 (2004), 34 -- 38.}






\end{thebibliography}
\end{document}